\def\versiondate{4 Aug. 2010}
\input math.macros
\input Ref.macros

\checkdefinedreferencetrue
\continuousfigurenumberingtrue
\theoremcountingtrue
\sectionnumberstrue
\forwardreferencetrue
\citationgenerationtrue
\nobracketcittrue
\hyperstrue
\initialeqmacro

\input\jobname.key
\bibsty{../../texstuff/myapalike}

\def\gh{G}
\def\vertex{{\ss V}}

\def\edges{{\ss E}}
\def\edge{{\ss E}}

\def\gp{\Gamma}
\def\gpe{\gamma}
\def\id{{\ss id}}  
\def\Aut{{\rm Aut}} 
\def\sqnm#1{\| #1\|_2^2}  
\def\PP{{\cal P}}   
\def\fh{\widehat f_B}  
\def\isoe{\Phi}  
\def\bfo{{\bf 1}}
\def\bfz{{\bf 0}}
\def\bde{\partial_\edges}
\def\F{\scr F}  

\ifproofmode \relax \else\head{{\it Europ. J. Combin.} {\bf 32} (2011), 1115--1125.}
{\versiondate}\fi 
\vglue20pt

\title{Perfect Matchings as IID Factors}
\title{on Non-Amenable Groups}

\author{Russell Lyons and Fedor Nazarov}

\abstract{We prove that in every bipartite
Cayley graph of every non-amenable group,
there is a perfect matching that is obtained as a factor of
independent uniform random variables. We also discuss expansion properties
of factors and improve the Hoffman spectral bound on independence number of
finite graphs.
}

\bottomIII{Primary 
37A50, 
22D40, 
22F10, 
60C05. 
Secondary
05C69, 
05C70.
}
{Independent sets, Hoffman, spectral bound, expansion, isoperimetric
constant.}
{Research partially supported by Microsoft Research and
NSF grants DMS-0406017 and DMS-0705518.}

\bsection{Introduction}{s.intro}

A \dfn{perfect matching} in a graph is a set of its edges that includes
each vertex exactly once.
An early result guaranteeing the existence of a perfect matching is due to
\refbmulti{Konig:finite,Konig:infinite}, who showed the sufficiency that
the graph be bipartite and regular of finite degree.
On the other hand, infinite graphs may come with a measurable structure and
one may wish for a measurable perfect matching.
That is, suppose that the bipartite graph has its two parts equal to $[0,
1]$ and $[2, 3]$, with edge set a Borel subset of $[0, 1] \times [2, 3]$.
If the graph is regular, must it have a Borel perfect matching?
\ref b.Lacz/ showed that the answer is no for 2-regular graphs.
\ref b.KNSS/ built on his example to show the same for any even
degree\ftnote{*}{One of us (R.L.) noted that their example was incorrect
for some degrees,
but their example was corrected by \ref b.ConKech/.}.
However, it is still open whether there is a measurable version of
K\"onig's theorem for 3-regular graphs.

A somewhat related notion is the following.
Suppose we are given a finitely generated group $\gp$ and a Cayley graph
$G$ of $\gp$.
In addition, we have independent uniform $[0, 1]$ random variables assigned
to each edge (or vertex). We call an instance of such an assignment a
\dfn{configuration}.
Note that $\gp$ acts on $G$ by automorphisms, whence it also acts on
the set of configurations, as well as on perfect matchings of $G$.
A random perfect matching of $G$ that is obtained as some measurable
function of the configuration and that commutes with the action
of $\gp$ is called a \dfn{$\gp$-factor} perfect matching of the random
variables.
Does one exist? In the case of the usual Cayley graph of $\Z$, the answer
is no since the only invariant measure on perfect matchings is not mixing,
yet every factor of independent random variables is mixing.
However, this may be the only exception.
\ref b.Timar:IM/ shows a positive answer for the usual Cayley graphs of
$\Z^d$ ($d > 1$).
Our main contribution is to prove that the answer is yes when $\gp$ is
non-amenable and $G$ is bipartite:

\procl t.matchingintro Let $G$ be a bipartite non-amenable simple Cayley
graph of a group $\gp$.
Then there is an $\gp$-factor
of independent uniform $[0, 1]$ random variables on $\gp$
that is a perfect matching of $G$ a.s.
\endprocl

\noindent
In fact, we prove a slight strengthening of this in \ref t.matching/.

One connection of the two above notions is due to
Kechris (personal communication, 2001; see \ref b.ConKech/).
He attempted to show that there is
no measurable version of K\"onig's theorem for 3-regular graphs by an
approach that would succeed if the 3-regular tree had no perfect matching
as a factor of IID. More precisely, note that the line graph $G'$ of the
3-regular tree is the
usual Cayley graph of $\Z_3 * \Z_3$ and that a perfect matching of the tree
corresponds to a set of vertices in $G'$ that has exactly one in each triangle.
His approach would succeed when there is no such set of vertices as a
$\Z_3 * \Z_3$-factor of independent random variables with values in $\{0,
1\}$.
This is equivalent to existence of a $\Z_3 * \Z_3$-factor from $[0,
1]^{\Z_3 * \Z_3}$ by a result of \ref b.Ball:factor/.
Our result is a factor that not only commutes with the action of the group,
but with all automorphisms of the Cayley graph, whence it shows that
Kechris's approach will not work, at least when sets of measure 0 are
ignored.

This is somewhat surprising, actually.
Consider again the case where $G$ is the 3-regular tree.
To obtain a perfect matching as a factor in $G$,
one must have a rule for each edge $e$ that decides whether $e$
belongs to the matching, depending on the configuration.
This rule must be the same for each edge (after action by an
automorphism); being measurable, it depends only on the configuration
within some distance $R$ of the edge, up to a small error.
The balls of radius $R$ about two neighboring edges have a
substantial symmetric difference, yet the rules must make consistent
decisions, so this appears very hard to do.
Indeed, if we wished to choose a set of vertices of $G$ as a factor with
the property that no two are adjacent, then there would be a bound to the
density of such a set that is strictly less than 1/2, even though $G$ is
bipartite.

As a matter of fact, there is considerable interest in finding such sets of
vertices with high density on regular trees.
The reason is this:
First, a set of vertices no two of which are adjacent is called
\dfn{independent}.
The \dfn{independence ratio} 
of a finite graph $\gh = (\vertex, \edge)$ is the maximum of $|K|/|\vertex|$
over all independent sets $K \subseteq \vertex$.
An open question is to determine the limiting independence ratio for
random regular graphs as the number of vertices tends to infinity.
The existence of this limit, but not its value, has recently been
established by \ref b.BGT/.
Any factor of IID random variables on a $d$-regular tree can be emulated on
finite $d$-regular graphs with large girth or, more generally, with
rare small cycles; this includes random regular graphs.
If the factor gives an independent set on the tree, then it will
give an independent set on the finite graphs of almost the same density.
Furthermore, the best lower bounds on the independence ratios on all
regular graphs of large girth are produced in this way
by factors on regular trees
(\ref B.LW/; \ref B.Hoppen:thesis/).
These match the best lower bounds on the independence ratios on random
regular graphs (see \ref b.Wormald:survey/), which were first obtained by other
techniques.
Furthermore, B.~Szegedy (personal communication, 2009) 
conjectures that the possible values for the
densities of independent sets in random $d$-regular graphs coincides with
the possible densities of independent sets as IID factors in $d$-regular
trees. (This is part of a much more general conjecture.)

In a wider context, factors are a fundamental object in the ergodic theory
of amenable groups. They are just beginning to be understood for
non-amenable groups: see \ref b.Bowen:markov/ for the case of free groups.

Finally, in the continuous context, Poisson point processes are the
analogue of IID random variables from the discrete setting and play a
corresponding role in the ergodic theory of continuous amenable groups.
There are several recent papers on factors of Poisson point processes that
give graphs, including perfect matchings: see \ref b.FLT:PT/, \ref
b.HP:treemat/, \ref b.Timar:pp/, and \ref b.HPPS:PM/.

In \ref s.spec/, we prove our result on perfect matchings.
This depends on an expansion property of factors in the non-amenable
setting. 
There, we also show how our theorem on perfect matchings extends to all
measure-preserving equivalence relations with expansive generating
graphings.
Some general remarks on expansion of factors are given in \ref
s.genl/.
Since matchings are independent sets in line graphs, we also discuss in
\ref s.hoff/ some
improvements in the classical Hoffman bound for independent sets.
We conclude with a few open questions in \ref s.question/.

\bsection{Perfect Matchings}{s.spec}

There are various equivalent definitions of non-amenability.
The simplest is due to \ref b.Folner/. To state it for a graph $\gh =
(\vertex,\edge)$,
define
$$
\isoe(\gh) := \inf \left\{ {|\bde K| \over |K|} \st
 \emptyset \ne K\subset \vertex \hbox{ is finite} \right\} \,.
$$
Here, $\bde K$ is the set of edges that join $K$ to its complement.
Then $G$ is \dfn{non-amenable} if $\isoe(G) > 0$.

We shall give a randomized
algorithm (that takes infinitely many steps) to produce
(a.s.) a perfect matching in a bipartite non-amenable Cayley graph.
To prove that it works, we shall need a lemma that exploits the expansion
property of non-amenability in the context of factors of IID.
Our proof of this expansion property depends on spectral information.

For a function $f : \gp \to \R$ and an element $\gpe \in \gp$, write
$R_\gpe f$ for the function $x \mapsto f(x \gpe)$.
The \dfn{right regular representation} of $\gp$ is the $\gp$ action on
$\ell^2(\gp)$ given by $\gpe \mapsto R_\gpe \restrict \ell^2(\gp)$.
A representation is called \dfn{subregular} if it is the restriction of the
regular representation to a $\gp$-invariant subspace.
The \dfn{trivial representation} is the action on $\R$ that fixes all
points.
Let $\mu$ be the usual product measure on $[0, 1]^\gp$ with each
coordinate getting Lebesgue measure.
We also have the representation $R$ on $L^2([0, 1]^\gp, \mu)$ given by
$(R_\gpe F)(\omega) := F(R_{\gpe^{-1}} \omega)$.
The following theorem has been known for some time; 
see Proposition 3.2 and Lemma 3.3 of \ref b.KT/
for a more general result.

\procl t.spectra
Let\/ $\gp$ be a countable group.
The representation $R$ of\/ $\gp$ on $L^2([0, 1]^\gp, \mu)$ is a sum of the
trivial representation, the regular representation, and subregular
representations.
\endprocl

(The basic idea is that if
$\{W_n\}$ is an orthonormal basis 
of $L^2\big([0, 1]\big)$ with $W_0 = \bfo$ and if $C_\gpe : [0, 1]^\gp
\to [0, 1]$ denotes the evaluation function at the coordinate $\gpe$,
then an orthonormal basis of $L^2\big([0, 1]^\gp, \mu^\gp\big)$ is the set
of all products $\prod_{\gpe \in \gp} W_{n(\gpe)} \circ C_\gpe$ with
$n(\gpe) = 0$ for all but finitely many $\gpe$.)

Fix a finite set $S \subset \gp$
that is closed under inverses and that generates $\gp$.
We are interested in the Cayley graph $\gh$ of $\gp$ with respect to $S$.
Let 
$$
\PP := |S|^{-1} \sum_{s \in S} R_s
$$
be the \dfn{transition operator}.
Then $\PP$ is self-adjoint and, thus, has real spectrum, whether it acts on
$\ell^2(\gp)$ or on 
$L^2([0, 1]^\gp, \mu)$.

The following is an immediate consequence of \ref t.spectra/.

\procl c.Pspectrum
The spectrum of\/ $\PP$ on $\ell^2(\gp)$
is the same as the spectrum of\/ $\PP$ on $\bfo^\perp$ in $L^2([0, 1]^\gp,
\mu)$.
\endprocl


Let $\rho$ be the spectral radius of $\PP$ on $\ell^2(\gp)$.
\ref b.Kesten:symm/ proved that $\rho < 1$ 
iff $\gp$ is nonamenable.
Let $X$ stand for $[0, 1]^{\gp}$.
Write $L^2_0(X, \mu)$ for the
orthocomplement of the constants in $L^2(X, \mu)$.

A measurable function $\phi : X \to \{0, 1\}^\gp$ or
$\phi : X \to \{0, 1\}^{\edge(G)}$ is called a \dfn 
{$\gp$-factor} if $R_\gpe \big(\phi (\omega)\big) = \phi\big(R_\gpe
\omega\big)$ for all $\gpe \in \gp$ and $\omega \in X$.
More generally, if $\gp'$ is a group of automorphisms of $G$ that commutes
with $\phi$, then $\phi$ is called a \dfn{$\gp'$-factor}.
The full group of automorphisms is denoted $\Aut(G)$.
To any factor with range $\{0, 1\}^\gp$, we associate the set
$$
B := \{ \omega \st \big(\phi(\omega)\big)(\id) = 1\}
\,,
$$
where $\id$ denotes the identity element of $\gp$.
Conversely, given any measurable $B \subseteq X$, there is an associated
$\gp$-factor defined by $\big(\phi(\omega)\big)(\gpe) := \II{R_\gpe \omega
\in B}$.
We think of the image $\phi(\omega)$ of a factor $\phi$
as subset of $\gp$, namely, those $\gpe \in \gp$ where
$\big(\phi(\omega)\big)(\gpe) = 1$ and also write $\gpe \in \phi(\omega)$
when $\big(\phi(\omega)\big)(\gpe) = 1$.
We sometimes omit parentheses and write $\phi\omega$ for $\phi(\omega)$.
We also think of
$$
b := \mu(B)
$$
as the density of the factor.
Write
$$
f_B := \I B - \mu(B) \bfo \in L^2_0(X, \mu)
\,.
$$
We have 
$$
\sqnm{f_B} = b (1-b)
\,.
$$

\procl l.nbrexpand
Let $G = (\gp, S)$ be a Cayley graph.
Let $\phi : (X, \mu) \to \{0, 1\}^\gp$
be a $\gp$-factor.
Define $\phi'\omega$ to consist of all the vertices that are adjacent to
some vertex in $\phi\omega$.
Then 
$$
b' := \P[\id \in \phi'\omega] \ge 
{1 \over \rho^2(1 - b) + b} \cdot b
\,.
\label e.indepnbrs
$$
\endprocl

\proof
Let $A := \{\omega \st \phi'\omega(\id) = 1\}$. 
Since $\I{A^c} \cdot \PP \I B = \bfz$,
we have 
$$
b
=
(\I B, \PP \bfo)
=
(\PP \I B, \bfo)
=
(\PP \I B, \I A)
\,.
$$
Therefore, 
$$
b^2
\le
\|\PP \I B\|^2 \|\I A\|^2
\,.
$$
Now $\|\I A\|^2 = b'$ and
$$
\|\PP \I B\|^2 
=
\|\PP (f_B + b \bfo)\|^2 
=
\|\PP f_B\|^2 + \|b \bfo\|^2 
\le
\rho^2 \|f_B\|^2 + b^2 
=
\rho^2 b(1-b) + b^2 
$$
since $f_B \perp \bfo$ and $\PP$ preserves $L^2_0(X, \mu)$.
Therefore, 
$$
{b' \over b}
\ge
{1 \over \rho^2(1 - b) + b}
\,.
\Qed
$$

%

We also need the following general tool 
(see, e.g., \ref b.BLPS:gip/), whose proof we include for the convenience
of the reader:

\proclaim The Mass-Transport Principle for Countable Groups.
Let\/ $\gp$ be a countable group. If $f : \gp \times \gp \to [0, \infty]$ is
diagonally invariant, then
$$
\sum_{x \in \gp} f(\id, x) = 
\sum_{x \in \gp} f(x, \id)
\,.
$$

\proof
Just note that $f(\id, x) = f(x^{-1}\id, x^{-1}x) = f(x^{-1}, \id)$ and
that summation of $f(x^{-1}, \id)$ over all $x^{-1}$ is the same as
$\sum_{x \in \gp} f(x, \id)$ since inversion is a bijection of $\gp$.
\Qed

In this context, we often use $f(x, y) = \E F(x, y, \omega)$, where $F$ is
defined on a probability space whose measure is $\gp$-invariant. If $F$ is
diagonally invariant, then so is $f$. We then call $F(x, y, \omega)$ the
\dfn{mass transported from $x$ to $y$}.

We are now ready to prove our main theorem.

\procl t.matching Let $G$ be a bipartite non-amenable simple Cayley graph.
Then there is an $\Aut(G)$-factor
of $\big([0, 1]^{\gp}, \mu\big)$
that is a perfect matching of $G$ a.s.
\endprocl

\proof
We shall construct the factor in infinitely many stages, each stage
consisting of infinitely many steps.
Since we can decompose a uniform [0, 1] random variable into an infinite
sequence of independent uniform [0, 1] random variables, we shall assume
that we are given such sequences at the start.
We shall also make use of a reverse operation: the \dfn{composition} of a
finite ordered list of numbers in [0, 1] is a number in [0, 1]. We choose
this composition map to be measurable and
so that given the length of the list of numbers, it is
an injection except on a countable set.
Each random variable will be used at most once.
We shall speak of the current random variables assigned to vertices, which we
throw away after use.

Suppose we have a (partial) matching.
Call a path \dfn{alternating} if its edges alternate between
belonging to the matching and not.
Following \ref b.Berge/, define an \dfn{(augmenting) chain} to be a simple
alternating path between unmatched vertices.
If we replace a chain by the same path, but with unmatched edges becoming
matched and matched edges becoming unmatched, so that all the vertices of
the path are now matched, we say that we \dfn{flip} the chain.

At the beginning of the first stage,
we have the empty matching and all edges are chains.
At the end of the $n$th stage, there will be no chains of length at most
$2n-1$, where length is measured by the number of edges.
Each step in the $n$th stage will be a repetition of the following
operation:
Assign the composition of the current random variables on the vertices to each
current chain of length at most $2n-1$.
If a chain has a larger composition than that of every other
chain that it intersects, then
flip that chain.

Note that once a vertex is matched in a given step, then it remains matched
after all subsequent steps.
Furthermore, each edge belongs to a finite number of chains of length at
most $2n-1$, whence it changes its status (between belonging to the
matching and not) at most finitely many times during the $n$th stage.
Finally, there is a lower bound (depending on $n$ and $|S|$) to the
conditional probability that a current chain is flipped, regardless of the
past, whence after infinitely many steps, there are a.s.\ no chains of
length at most $2n-1$.

In order to define the factor as a limit of the stages, we must prove that
a.s.\ no edge changes its status infinitely many times.

Let $\phi_n$ denote the factor matching at the end of the $n$th stage.
Fix $n$ and
define $\Seq{A_k}$ recursively as follows.
Let $A_0 = A_0(\omega)$ denote the unmatched vertices in $\phi_n \omega$.
If $k$ is even, then let $A_{k+1}$ be the set of vertices that have a
neighbor in $A_k$.
If $k$ is odd, then let $A_{k+1}$ be the set of vertices $x$ such that for
some $y \in A_k$, the edge $[x, y]$ is matched in $\phi_n \omega$.

We claim that for every $k \ge 1$ and every $x \in A_k$, there is a simple
alternating path from some $x_0 \in A_0$ to $x$ of length at most $k$.
Clearly, there is some alternating path $P_x$ from some $x_0 \in A_0$ to
$x$ of length at most $k$.
Since $G$ is bipartite, each edge of $P_x$ that leads from a vertex at odd
distance from $x_0$ to a vertex at even distance from $x_0$ is a
matched edge, whence the shortest path from $x_0$ to $x$ contained in $P_x$
is simple and alternating. 

There are two consequences of this that we use: The first is that if $x \in
A_k$ is unmatched and $k \ge 1$, then there is a chain of length at most $k$.
The second is that
if for some even $k$, the set $A_k$ is not independent, then
there is a chain of length at most $2k+1$.
Indeed, suppose that $x, y \in A_k$ are neighbors.
By the above, there is some simple alternating path $P_x$ from some
$x_0 \in A_0$ to $x$ of even length at most $k$
and a simple alternating path $P_y$ from some $y_0 \in A_0$
to $y$ of even length at most $k$.
Since the concatenation $P$ of $P_x$ followed by the edge $(x, y)$ and
then finally the reverse of $P_y$ is a path of odd length from
$x_0$ to $y_0$, it follows that the distance between $x_0$ and $y_0$ is
odd. In particular, $x_0 \ne y_0$.
If $P_x$ and $P_y$ are disjoint, then since the last edge of each of these
paths lies in the matching, the path $P$ is a chain of length at most
$2k+1$, as desired.
In case $P_x$ and $P_y$ are not disjoint, then their union contains a
simple path $Q$ from $x_0$ to $y_0$.
Since the length of $Q$ is odd, it is easy to see that $Q$ is alternating
as well.

By the first consequence, when $k < n$ is odd, there is a unique edge in
the matching from each $x \in A_k$ to some $y \in A_{k+1}$.
Let $x$ send mass 1 to $y$ in this situation.
Then by the Mass-Transport Principle,
$\mu[\id \in A_{k+1}] = \mu[\id \in A_k]$ for all odd $k < n$, where
$\mu[\id \in A_k]$ means $\mu\big(\{\omega \st \id \in A_k(\omega)\}\big)$.
By the second consequence,
for all even $k < n$, the set $A_k$ is independent,
which implies (for example, by \ref l.nbrexpand/)
that $\mu[\id \in A_k] \le 1/2$.
By \ref l.nbrexpand/, $\mu[\id \in A_{k+1}] \ge c \mu[\id \in A_k]$ for all even
$k < n$, where $c := 2/(1+\rho^2)$. Note that
$c > 1$ because $G$ is non-amenable.
(If $G$ is a tree, then instead of using \ref l.nbrexpand/, one could
deduce this expansion inequality for $\mu[\id \in A_k]$ by using
the fact that regular trees are limits of finite
bipartite expander graphs, in the sense that the proportion tends to 1
of vertices in those finite graphs with a large neighborhood the same as in
the tree.)
Since 
$\mu[\id \in A_{2k-1}] \ge c^k \mu[\id \in A_0]$ for $2k-1 \le n$, it
follows that $\mu[\id \in A_0] \le a_n := c^{-\flr{(n+1)/2}}$.

Now let each endpoint of a chain that is flipped 
send mass 1 to each vertex in its flipped chain.
Then the expected mass sent by the identity is at most $\sum_n 2n
a_{n-1} < \infty$.
Each vertex receives mass equal to twice the number of times a neighboring
edge changes its status.
By the Mass-Transport Principle, the expected number of times an edge
changes its status is finite.
This proves that the limit of $\phi_n$ exists a.s.\ and that a.s.\
all vertices are matched at the end.
\Qed

\procl r.edgefactor
The same result holds for factors from
$\big([0, 1]^{\edge(G)}, \nu\big)$, where $\nu$ is product measure. Indeed 
$\big([0, 1]^{\gp}, \mu\big)$ is itself a factor of
$\big([0, 1]^{\edge(G)}, \nu\big)$.  
To see this, given $\omega \in 
[0, 1]^{\edge(G)}$, define $\xi \in   
[0, 1]^{\gp}$ by $\xi(x) := \sum_{e \ni x} \omega(e)$
(mod 1).
It is clear that each $\xi(x)$ is uniform on $[0, 1]$ when $\omega$ has law
$\nu$.
To prove that $\xi(x_1), \ldots, \xi(x_n)$ are independent, we proceed by
induction.
Because $\gh$ is infinite, we may assume that $x_n$ belongs to an edge $e$
whose other endpoint is not among $x_1, \ldots, x_{n-1}$.
Since $\omega(e)$ is therefore independent of $\xi(x_i)$ for $i < n$, it
follows that $\xi(x_n)$ is independent of $\xi(x_i)$ for $i < n$.
\endprocl

\procl r.mper
Let $(X, \F, \mu)$ be a probability space and $E \in \F \times \F$ a
symmetric measurable subset of $X \times X$. Let $G := (X, E)$ be the graph
associated to $E$. Assume that all the connected components of $G$ are
bipartite and denumerable. Write $[x] \subset X$ for the vertices in the
connected component of $x \in X$.
Suppose that $(X, \F, \mu, G)$ is measure-preserving, 
meaning that $\mu_{\rm L} = \mu_{\rm R}$, where
$$
\int_{X^2} f(x, y) \,d\mu_{\rm L}(x, y)
:= \int_{x \in X} \sum_{y \in [x]} f(x, y) \,d\mu(x)
$$
and
$$
\int_{X^2} f(x, y) \,d\mu_{\rm R}(x, y)
:=
\int_{x \in X} \sum_{y \in [x]} f(y, x) \,d\mu(x)
$$
for all measurable $f : X^2 \to [0, \infty]$.
Suppose in addition that $G$ is expansive,
meaning that there exists $c > 1$ such that for every measurable $A \subset
X$ with $\mu(A) \le 1/2$, we have $\mu(A') \ge c \mu(A)$, where $A'$
consists of the $G$-neighbors of the points in $A$. Then there is a
$\mu_{\rm L}$-measurable perfect matching in $G$. 
The proof is the same as that of \ref t.matching/, except that the first
short part is replaced by a (similar) general argument of \ref b.ElekLipp/,
Proposition 1.1,
which shows that there is a sequence of factors $\phi_n$ that have the
property that there is no chain of length at most $2n-1$ in $\phi_n$ and
such that the set of matched vertices is increasing in $n$.
\endprocl

\comment{
As a consequence of \ref t.matching/ and \ref p.independent/, we have that
for the usual Cayley graph of $\gp := \Z_3 * \Z_3$,
there is an eigenfunction of $\PP$ on $L^2\big([0, 1]^\gp, \mu\big)$
for the eigenvalue $-\rho_- = -1/2$ that takes only two values, namely,
$f_B$ for the perfect matching factor $B$.
}

\bsection{Factor Expansion and Spectral Radius}{s.genl}

There is a general relationship between factors of measure-preserving
actions and an associated spectral radius.
It is quite analogous to expansion properties of finite graphs.

Let $\gp$ be a group acting by measure-preserving transformations on a
probability space $(X, \mu)$.
We also write integration with respect to $\mu$ as $\E$.
Fix a finite $S \subset \gp$,
closed under inverses and generating $\gp$.
Let $\rho$ be the spectral radius of $\PP$ on $L^2_0(X, \mu)$.
In fact, for more precision, we shall use the bottom, $-\rho_-$, 
and the top, $\rho_+$, of the spectrum. We have $\rho = \max (\rho_-,
\rho_+)$ and
$$
-\rho_- \le (\PP f, f) \le \rho_+ 
\label e.specint
$$
for all $f \in L^2_0(X, \mu)$ with $\|f\|_2 = 1$.


Define the \dfn{expansion constant} of the action with respect to $S$ by 
$$
\isoe(\gp, S, X, \mu)
:=
\inf \Big\{ {1 \over |S| b (1-b)} \int \sum_{s \in S} \I{B \cap s B^c}
\,d\mu \st B \subset X,\, 0 < \mu B < 1 \Big\}
\,.
$$

The following inequalities relating the expansion constant and the spectral
radius are analogous to those on finite graphs, so we restrict our proofs
to the essential steps. See, e.g., \ref b.LPW:book/, Theorem 13.14, for
more details on finite graphs.

\procl t.genl
Let $\gp$ be a group acting by measure-preserving transformations on a
probability space $(X, \mu)$.
Write $\isoe := \isoe(\gp, S, X, \mu)$.
Then
$$
\isoe^2 / 8
\le
1 - \sqrt{1 - (\isoe/2)^2}
\le
1 - \rho
\le
1 - \rho_+
\le
\isoe
\,.
$$
\endprocl

There is never expansion for amenable groups, that is, for all
actions of an amenable group, the spectral radius is equal to 1 by a
theorem of \ref b.OrnWei:Rohlin/.
Some such expansion property holds for every ergodic invariant percolation
only on Kazhdan groups.
In fact, the very definition of Kazhdan's property (T) is easily seen to be
equivalent to every action having spectral radius less than 1.
As we saw via \ref c.Pspectrum/,
the spectral radius of Bernoulli actions of non-amenable groups
is strictly less than 1.

For some purposes, a notion for an action weaker than expansion is
interesting, namely, the non-existence of almost invariant sets. This means
that if $\gp$ acts on $(X, \mu)$ and
$B_n \subset X$ are measurable sets with $\mu(B_n \triangle \gpe B_n) \to
0$ for all $\gpe \in \gp$, then $\mu(B_n)\big(1 - \mu(B_n)\big) \to 0$. See
Appendix A of \ref b.HjKech:rigid/ for a discussion of this and related
matters.

We have 
$$
\isoe
=
\inf {1 \over b (1-b)} (\I B, \PP \I{B^c})
=
1 - \sup {(\PP f_B, f_B) \over b(1-b)}
\,,
$$
which proves that $\isoe \ge 1 - \rho_+$.
Since $(\I B, \PP \I{B^c}) = (\PP \I B, \I{B^c})$,
we also have the alternative expression 
$$
\isoe = \inf \Big\{ {1 \over 2 |S| b (1-b)} \E \sum_{s \in S} 
(\I{B \cap s B^c} + \I{B^c \cap s B})
\st B \subset X,\, 0 < \mu B < 1 \Big\}
\,.
\label e.alt
$$

\procl l.levelsets
With notation as in \ref t.genl/,
if $f \in L^2(X, \mu)$ satisfies $f \ge 0$ a.s., then
$$
2 \mu[f = 0] \isoe \int f \,d\mu
\le
{1 \over |S|} \int \sum_{s \in S} |f - s f| \,d\mu
\,.
$$
\endprocl

\proof
For $t > 0$, let $B_t := f^{-1}(t, \infty)$.
Put $\alpha_f := \mu[f = 0]$.
Then by \ref e.alt/, we have 
$$
2 \isoe \mu(B_t) \alpha_f |S|
\le
2 \isoe \mu(B_t) \mu(B_t^c) |S|
\le
\E \sum_s \big(\II{f > t \ge s f} + \II{s f > t \ge f}\big) 
\,.
$$
Integrating on $t \in (0, \infty)$ with respect to Lebesgue measure gives 
$$
2 \isoe \alpha_f |S| \E f
\le
\E \sum_s \big( \max\{f - s f, 0\} + \max\{ s f - f, 0\}\big)
=
\E \sum_s |f - s f|
\,.
\Qed
$$

\proofof t.genl
We have already proved the fourth inequality. The first inequality
is elementary. To prove the second inequality, consider $f_0 \in L^2_0$
such that $\|f_0\| = 1$.
Define $\lambda := (\PP f_0, f_0)$.
Without loss of generality, we may assume that $\mu[f_0 > 0] \le 1/2$.
Define $f := \max\{f_0, 0\}$.
Then checking cases shows that $(I - \PP) f \le (1 - \lambda) f$, whence
$\big((I - \PP)f, f\big) \le (1 - \lambda) \|f\|^2$.
Define 
$$
\beta := \E \sum_s |f - s f|^2/(2|S|)
=
\big((I - \PP)f, f\big)
\,.
$$
Now by the lemma, since $\alpha_f \ge 1/2$, we have 
$$
\|f\|_2^4
\le
\isoe^{-2} \big(\E \sum_s |f^2 - s f^2|/|S|\big)^2
\le
2 \isoe^{-2} \beta \E \sum_s |f + s f|^2/|S| 
=
2 \isoe^{-2} \beta (4 \|f\|_2^2 - 2 \beta)
\,.
$$
Therefore, 
$$
\lambda^2
\le
\big(1 - \beta/\|f\|_2^2\big)^2
\le
1 - (\isoe/2)^2
\,.
$$
Now take the supremum of $|\lambda|$ over $f_0$.
\Qed

Let $G := (\gp, S)$ be the right Cayley graph of $\gp$ corresponding to the
generating set $S$.
%
%
%
When the factor is an independent set in $G$, we can bound its density as
follows.
It is analogous to the Hoffman bound (\ref B.Lov:shannon/)
for the independence number of a finite graph.

\procl p.independent
Suppose that $\phi : (X, \mu) \to \{0, 1\}^\gp$ is a $\gp$-factor with
the property that if $(\phi\omega)(\id) = 1$, then $(\phi\omega)(s) = 0$
for all $s \in S$.
Then 
$$
b \le \rho_-/(1+\rho_-)
\,,
\label e.indepbd
$$
with equality iff $\PP f_B = -\rho_- f_B$.
\endprocl

\proof
We have $(\PP \I B, \I B) = 0$,
which is the same as 
$$
(\PP f_B, f_B) = - b^2
\,.
\label e.inp
$$
We deduce from this that $b^2 \le b(1-b)\rho_-$, which gives the inequality.
Furthermore, if equality holds, then
$$
-\rho_- = {(\PP f_B, f_B) \over \sqnm{f_B}}
\,.
$$
By \ref e.specint/, it follows that $\PP f_B = -\rho_- f_B$.
Conversely, if $\PP f_B = -\rho_- f_B$, then it easily follows that
equality holds in \ref e.indepbd/.
\Qed

%
%



Other inequalities known for finite graphs can be proved as well. We
illustrate with two well-known examples (see, e.g., \ref
b.AlonSp/, Theorem 9.2.4 and Corollary 9.2.5).

\procl p.variance
Let $G = (\gp, S)$ be a Cayley graph.
Let $\phi : (X, \mu) \to \{0, 1\}^\gp$ be a $\gp$-factor.
Then 
$$
\Eleft{\left[{1 \over |S|} \sum_{s \in S} \I{s B} - b\right]^2}
\le
\rho^2 b (1-b)
\,.
$$
\endprocl

\proof
This is the same as $\|\PP f_B\|^2 \le \rho^2 \|f_B\|^2.$
\Qed

\procl c.expmix
Let $G = (\gp, S)$ be a Cayley graph.
Let $\phi : (X, \mu) \to \big(\{0, 1\}^2\big)^\gp$ be a $\gp$-factor.
Define $B_i := \{\omega \st \big(\phi\omega(\id)\big)_i = 1\}$ for $i = 1, 2$.
Put $b_i := \mu(B_i)$.
Then 
$$
\left|\Eleft{{1 \over |S|} \sum_{s \in S} \I{B_1 \cap s B_2}} - b_1 b_2\right|
\le
\rho \sqrt{b_1 b_2 (1-b_1) (1 - b_2)}
\,.
$$
\endprocl

\proof
The left-hand side equals $|(f_{B_1}, \PP f_{B_2})|$, whence it is at most
$\|f_{B_1}\| \cdot \|\PP f_{B_2}\|$.
Multiplying this by the same inequality with $B_1$ and $B_2$ reversed and
using \ref p.variance/ gives the result.
\Qed

\bsection{Improving the Hoffman Bound}{s.hoff}

Here we discuss briefly how to improve \ref p.independent/. Our results
apply to factors as well as to arbitrary regular finite graphs.
One improvement holds only when $\rho_- > 1 - 1/|S|$; the other holds
when $|S|\rho_- \notin \Z$. 
In both cases, we give only sketches since we have no especially
interesting applications to present.
However, since the Hoffman bound has not been improved since its discovery,
it seems worthwhile to explain our improvements.

There are various ways to improve the proof of \ref p.independent/. We give
just one.
Given the factor $\phi$ such that $\phi\omega$ is a.s.\ an independent set,
define $N(\omega) := |S \cap \phi\omega|$.
Write $d := |S|$ and $p := \P[N = d]$.
We may assume that $\phi\omega$ is a.s.\ a maximal independent set, i.e.,
every vertex not in $\phi\omega$ has a neighbor in $\phi\omega$.
Consider the function 
$$
f(\omega) := \cases{\phantom{-}1 &if $N(\omega) = 0$,\cr
                    -a &if $1 \le N(\omega) \le d- 1$,\cr 
                    -A &if $N(\omega) = d$.}
$$
We choose the values of $a$ and $A$ so that $f \perp \bfo$.
Then 
$$
d\cdot \PP f(\omega) =
\cases{-a d + (a - A) |\{s \in S \st N(s\omega) = d\}| &if $N(\omega) = 0$,\cr
       (1+a) N(\omega) - a d &if $N(\omega) \ge 1$.}
$$
Using the facts that $\E[N] = d b$, 
$$
\E[N \st 1 \le N \le d - 1] = \E[N] - d p = d (b - p)
\,,
$$
and $\Ebig{|\{s \in S \st N(\omega) = 0,\, N(s\omega) = d\}|} =
\Ebig{|\{s \in S \st N(s\omega) = d\}|} = d p$, one can calculate that 
$$
(\PP f, f) 
=
(1 + a)^2(1 - 2b) - 1
\,.
$$
Also, 
$$
(f, f) = b + a^2(1-b-p) + A^2 p
=
b + a^2(1-b-p) + [b - a(1 - b - p)]^2/p
\,.
$$
Since $(\PP f, f) \ge -\rho_- (f, f)$ for all $a$ and this inequality is
quadratic in $a$, it follows that its discriminant is non-positive: 
$$
0 \ge b \rho_-[b(1+\rho_-) - \rho_-] + p[1 - b(1+\rho_-)(2-\rho_-)]
\,.
$$
Since $b \le \rho_-/(1+\rho_-) \le 1/(1+\rho_-)(2-\rho_-)$, the same
inequality holds when we substitute a lower bound for $p$.
Now $\E[N \st 1 \le N \le d-1] \le (d-1)(1-b-p)$, which yields $p \ge
(2d-1)d - d + 1$, whence 
$$
b \le {d - 1 \over (1+\rho_-)[d(2-\rho_-) - 1]}
\,.
\label e.hoff1
$$
As we said, \ref e.hoff1/ improves \ref p.independent/ only for $\rho_- > 1
- 1/d$. In fact, it is impossible to improve \ref p.independent/ in all
cases when $\rho_- = 1 - 1/d$ because there are cases when equality holds.

Our second improvement is as follows.
We have 
$$
\PP \I B = {q \over d} \I {B^c}
\label e.defq
$$
for some integer-valued function $q$.
Now
$$
\Ebig{q/|S| \bigm| B^c} = b/(1-b)
\,.
\label e.condexp
$$
Write $\fh := f_B/\|f_B\|$.
If $\nu$ denotes the spectral measure for $\fh$ with respect to $\PP$, 
then 
$$\eqaln{
\|\PP \fh + \rho_- \fh\|^2
&=
\int_{-\rho_-}^{\rho_+} (\lambda + \rho_-)^2 \,d\nu(\lambda)
\cr&\le
(\rho_- + \rho_+)
\int_{-\rho_-}^{\rho_+} (\lambda + \rho_-) \,d\nu(\lambda)
\cr&=
(\rho_- + \rho_+)
[- {b \over 1 - b} + \rho_-]
\,.
}$$
On the other hand, by \ref e.defq/, one can calculate that
$$
b(1-b)\|\PP \fh + \rho_- \fh\|^2
=
\Var(q/d \mid B^c) (1-b) + [\rho_-(1-b) - b]^2 b/(1-b)
\,,
$$
whence 
$$
(\rho_- + \rho_+) b [\rho_-(1-b) - b]
\ge
\Var(q/d \mid B^c) (1-b) + [\rho_-(1-b) - b]^2 b/(1-b)
\,,
$$
which simplifies to 
$$
\Var\big(q/|S| \bigm| B^c\big) \le
b \left[\rho_- - {b \over (1-b)}\right]\left[\rho_+ + {b \over (1-b)}\right]
\,.
\label e.var
$$

Now if $m := b/(1-b) \in (k/d, (k+1)/d)$, then
the smallest $\Var(q/d \mid B^c)$ can be is when $q$
takes only the values $k$ and $k+1$ on $B^c$.
Note that if $\rho_- < (k+1)/d$, then $m < (k+1)/d$.
This gives that either $m \le k/d$ or 
$$
\Var(q/d \mid B^c)
\ge
-{k + k^2 \over d^2} + {1 + 2 k \over d} m - m^2
\,.
$$
Combining this with \ref e.var/ and, for simplicity, using $\rho \le 1$, we
obtain 
$$
m
\le
{k^2+k \over d(1+2k) - d^2 \rho_-}
\,,
$$
which is the same as 
$$
b
\le
{k^2+k \over k^2 + (1+2d)k + d - d^2 \rho_-}
\,.
$$

\procl r.var
Expanding the inequality of
\ref p.variance/ gives the same inequality \ref e.var/, but with
$\rho$ in place of $\rho_-$, which can be significantly worse.
\endprocl

\comment{
\procl p.irrat
If $\phi : (X, \mu) \to \{0, 1\}^\gp$ is a $\gp$-factor that is an
independent set a.s.\ in the Cayley graph of $(\gp, S)$,
then there is a measurable function $q : X \to \big[0,
|S|\big] \cap \Z$ such that 
$$
\Ebig{q/|S| \bigm| B^c} = b/(1-b)
\label e.condexp
$$
and
$$
\Var\big(q/|S| \bigm| B^c\big) \le
b \left[\rho_- - {b \over (1-b)}\right]\left[\rho_+ + {b \over (1-b)}\right]
\,.
\label e.var
$$
\endprocl


\proof
Write $d := |S|$ for the degree.
We have 
$$
\PP \I B = {q \over d} \I {B^c}
\label e.defq
$$
for some integer-valued function $q$.
Now 
$$
b = \int \PP \I B \,d\mu
=
\int_{B^c} {q \over d} \,d\mu
\,,
$$
whence 
\ref e.condexp/ holds.
Let $E_\lambda$ be the spectral measure for $\PP$: 
$$
\PP = \int_{-\rho_-}^{\rho_+} \lambda \,dE_\lambda
\,.
$$
Write $\fh := f_B/\|f_B\|$.
Then 
$$
\int_{-\rho_-}^{\rho_+} (\lambda + \rho_-) \,dE_\lambda \fh
=
\PP \fh + \rho_- \fh
\,,
$$
whence
$$\eqaln{
\|\PP \fh + \rho_- \fh\|^2
&=
\|\int_{-\rho_-}^{\rho_+} (\lambda + \rho_-) \,dE_\lambda \fh\|^2
\cr&=
\int_{-\rho_-}^{\rho_+} (\lambda + \rho_-)^2 \,d(E_\lambda \fh, \fh)
\cr&\le
(\rho_- + \rho_+)
\int_{-\rho_-}^{\rho_+} (\lambda + \rho_-) \,d(E_\lambda \fh, \fh)
\cr&=
(\rho_- + \rho_+)
[(\PP \fh, \fh) + \rho_-]
=
(\rho_- + \rho_+)
[- {b \over 1 - b} + \rho_-]
}$$
by \ref e.inp/.
On the other hand, by \ref e.defq/, we have 
$$\eqaln{
b(1-b)\|\PP \fh + \rho_- \fh\|^2
&=
\|(q/d) \I{B^c} - b + \rho_- \I B - \rho_- b\|^2
\cr&=
\E[((q/d) - (\rho_-+1)b)^2 ; B^c] + [\rho_-(1-b) - b]^2 b
\cr&=
[\Var(q/d \mid B^c) + (b/(1-b) - (\rho_-+1)b)^2] (1-b) 
+ [\rho_-(1-b) - b]^2 b
\cr&=
\Var(q/d \mid B^c) (1-b) + [\rho_-(1-b) - b]^2 b/(1-b)
\,,
}$$
whence 
$$
(\rho_- + \rho_+) b [\rho_-(1-b) - b]
\ge
\Var(q/d \mid B^c) (1-b) + [\rho_-(1-b) - b]^2 b/(1-b)
\,,
$$
which simplifies to \ref e.var/.
\Qed

Since $|q/d - b/(1-b)|$ is bounded away from 0 for $b$ close to
$\rho_-/(1+\rho_-)$, the inequality \ref e.var/ is impossible
if $b$ is too close to the bound in \ref p.independent/.

For example, consider
the Cayley graph of\/ $\Z * \Z_4$ with respect
to the usual generating set, so that it has degree 4.
A formula for the spectral radius $\rho$ was given by \ref b.Paschke/;
since the graph is bipartite, computation yields 
$\rho_+ = \rho_- = 0.89237^+$.
Now if $m := b/(1-b) \in (3/4, 1)$, then
the smallest $\Var(q/d \mid B^c)$ can be is when $q$
takes only the values 3 and 4 on $B^c$.
Similarly for other possible values of the mean.
This gives that
$$
\Var(q/d \mid B^c)
\ge
\max[m/4, (6 m - 1)/8, (10 m - 3)/8, (7 m - 3)/4] - m^2
\,.
$$
Following the arithmetic gives 
$$
b < 0.4663 < \rho/(1+\rho) = 0.47156^+
\,.
$$
\msnote{Choose a larger degree example to do better than our new spectral
bound.}

\procl r.var
Expanding the inequality of
\ref p.variance/ gives the same inequality \ref e.var/, but with
$\rho$ in place of $\rho_-$, which can be worse.
\endprocl
}

\bsection{Open Questions}{s.question}

It is interesting to consider the chromatic number with respect to
invariant processes under increasing restrictions:
For example, a regular tree has chromatic number 2, and there is an invariant
random proper 2-coloring, which is ergodic.
However, there is no such mixing 2-coloring, but there is a mixing
3-coloring.
What is the minimum number of colors for a proper coloring that is an
IID factor? For large degree $d$, it is at least $d/(2 \log d)$ since 
\ref b.FrLu/ proved that for large degree $d$, 
the independence ratio
for large random $d$-regular graphs is asymptotic to
$2 \log d/d$.
The minimum number of colors as an IID factor on a Cayley graph of degree
$d$ is at most $d+1$, as shown by Schramm (personal communication, 1997).
This was shown more generally to hold for any factor by \ref b.KST:color/.

Related to this is a simpler question due to Lyons and Schramm in 1997
(unpublished):
If $\gh$ is a Cayley graph of
chromatic number $\chi$, then is there a random invariant $\chi$-coloring?
It is easy to show a positive answer when $\gh$ is amenable.
\ref b.ConKech/ prove some general results on invariant coloring, as
well as a version of our \ref p.independent/, discovered independently.
For various results on coloring Poisson-Voronoi tessellations
by factors, see \ref b.ABGMP/ and \ref b.Timar:vor/.



\medbreak
\noindent {\bf Acknowledgements.}\enspace We are very grateful to David
Aldous, Yuval Peres, and Alexander Kechris for asking questions that
inspired this work and to Benjy Weiss for describing a version of \ref
t.spectra/.

\def\noop#1{\relax}
\input \jobname.bbl

\filbreak
\begingroup
\eightpoint\sc
\parindent=0pt\baselineskip=10pt

Department of Mathematics,
831 E. 3rd St.,
Indiana University,
Bloomington, IN 47405-5701
\emailwww{rdlyons@indiana.edu}
{http://mypage.iu.edu/\string~rdlyons/}

and

Department of Mathematics,
University of Wisconsin,
480 Lincoln Drive,
Madison, WI 53706 
\emailwww{nazarov@math.wisc.edu}
{http://www.math.wisc.edu/\string~nazarov/}


\endgroup


\def\cprime{$'$} \def\cprime{$'$} \def\cprime{$'$}
\def\temp{\let\linkit=\linkyear \apaliketrue}
\temp
\ifcitationgeneration\immediate\write\labelfile{\sanitize\temp}\fi
\def\startreferences{
 \vskip0pt plus.3\vsize \penalty -150 \vskip0pt
 plus-.3\vsize \bigskip\bigskip \vskip \parskip
 \begingroup\baselineskip=12pt\frenchspacing
 \bibliographytitle
 \vskip12pt\parindent=0pt
 \def\and{{\rm and}}
 \def\em{\it}
 \def\newblock{\hskip .11em plus.33em minus.07em}
 \def\bibauthor##1{{\sc ##1}}
 \def\bibitem[##1]##2
 {\htmlanchor{##2}{}\RefLabel{##2}[##1]\hangindent=.8cm\hangafter=1}
 }
\def\endreferences{\bigskip\bigskip\endgroup}
\ifundefined{bibstylemodification}\relax\else\bibstylemodification\fi
\startreferences

\bibitem[Alon and Spencer (2008)]{MR2437651}
\bibauthor{Alon, N. \and{} Spencer, J.H.} (2008).
\newblock {\em The Probabilistic Method}.
\newblock Wiley-Interscience Series in Discrete Mathematics and Optimization.
  John Wiley \& Sons Inc., Hoboken, NJ, third edition.
\newblock With an appendix on the life and work of Paul Erd{\H{o}}s.

\bibitem[Angel, Benjamini, Gurel-Gurevich, Meyerovitch, and Peled
  (2009)]{ABGMP}
\bibauthor{Angel, O., Benjamini, I., Gurel-Gurevich, O., Meyerovitch, T.,
  \and{} Peled, R.} (2009).
\newblock Stationary map coloring.
\newblock Preprint, \arXiv{0905.2563}.

\bibitem[Ball (2005)]{MR2142942}
\bibauthor{Ball, K.} (2005).
\newblock Factors of independent and identically distributed processes with
  non-amenable group actions.
\newblock {\em Ergodic Theory Dynam. Systems} {\bf 25}, 711--730.

\bibitem[Bayati, Gamarnik, and Tetali (2010)]{BGT}
\bibauthor{Bayati, M., Gamarnik, D., \and{} Tetali, P.} (2010).
\newblock Combinatorial approach to the interpolation method and scaling limits
  in sparse random graphs.
\newblock {\em Proceedings of the 42nd {ACM} symposium on Theory of Computing},
  105--114.
\newblock Cambridge, Massachusetts, USA.

\bibitem[Benjamini, Lyons, Peres, and Schramm (1999)]{MR99m:60149}
\bibauthor{Benjamini, I., Lyons, R., Peres, Y., \and{} Schramm, O.} (1999).
\newblock Group-invariant percolation on graphs.
\newblock {\em Geom. Funct. Anal.} {\bf 9}, 29--66.

\bibitem[Berge (1957)]{MR0094811}
\bibauthor{Berge, C.} (1957).
\newblock Two theorems in graph theory.
\newblock {\em Proc. Nat. Acad. Sci. U.S.A.} {\bf 43}, 842--844.

\bibitem[Bowen (2010)]{Bowen:markov}
\bibauthor{Bowen, L.} (2010).
\newblock Non-abelian free group actions: {M}arkov processes, the
  {A}bramov-{R}ohlin formula and {Y}uzvinskii's formula.
\newblock {\em Ergodic Theory Dynamical Systems}.
\newblock Available on CJO 13 Oct 2009, doi:10.1017/S0143385709000844.

\bibitem[Conley and Kechris (2009)]{ConKech}
\bibauthor{Conley, C.T. \and{} Kechris, A.S.} (2009).
\newblock Measurable chromatic and independence numbers for ergodic graphs and
  group actions.
\newblock Preprint.

\bibitem[Elek and Lippner (2010)]{ElekLipp}
\bibauthor{Elek, G. \and{} Lippner, G.} (2010).
\newblock Borel oracles. {A}n analytical approach to constant-time algorithms.
\newblock {\em Proc. Amer. Math. Soc.} {\bf 138}, 2939--2947.

\bibitem[Ferrari, Landim, and Thorisson (2004)]{MR2044812}
\bibauthor{Ferrari, P.A., Landim, C., \and{} Thorisson, H.} (2004).
\newblock Poisson trees, succession lines and coalescing random walks.
\newblock {\em Ann. Inst. H. Poincar\'e Probab. Statist.} {\bf 40}, 141--152.

\bibitem[F{\o}lner (1955)]{MR18:51f}
\bibauthor{F{\o}lner, E.} (1955).
\newblock On groups with full {B}anach mean value.
\newblock {\em Math. Scand.} {\bf 3}, 243--254.

\bibitem[Frieze and {\L}uczak (1992)]{MR1142268}
\bibauthor{Frieze, A.M. \and{} {\L}uczak, T.} (1992).
\newblock On the independence and chromatic numbers of random regular graphs.
\newblock {\em J. Combin. Theory Ser. B} {\bf 54}, 123--132.

\bibitem[Hjorth and Kechris (2005)]{MR2155451}
\bibauthor{Hjorth, G. \and{} Kechris, A.S.} (2005).
\newblock Rigidity theorems for actions of product groups and countable {B}orel
  equivalence relations.
\newblock {\em Mem. Amer. Math. Soc.} {\bf 177}, viii+109.

\bibitem[Holroyd, Pemantle, Peres, and Schramm (2009)]{HPPS:PM}
\bibauthor{Holroyd, A.E., Pemantle, R., Peres, Y., \and{} Schramm, O.} (2009).
\newblock Poisson matching.
\newblock {\em Ann. Inst. H. Poincar\'e Probab. Statist.} {\bf 45}, 266--287.

\bibitem[Holroyd and Peres (2003)]{MR2004b:60127}
\bibauthor{Holroyd, A.E. \and{} Peres, Y.} (2003).
\newblock Trees and matchings from point processes.
\newblock {\em Electron. Comm. Probab.} {\bf 8}, 17--27 (electronic).

\bibitem[Hoppen (2008)]{Hoppen:thesis}
\bibauthor{Hoppen, C.} (2008).
\newblock {\em Properties of Graphs with Large Girth}.
\newblock Ph.D. thesis, University of Waterloo.

\bibitem[Kechris, Solecki, and Todorcevic (1999)]{MR1667145}
\bibauthor{Kechris, A.S., Solecki, S., \and{} Todorcevic, S.} (1999).
\newblock Borel chromatic numbers.
\newblock {\em Adv. Math.} {\bf 141}, 1--44.

\bibitem[Kechris and Tsankov (2008)]{MR2358510}
\bibauthor{Kechris, A.S. \and{} Tsankov, T.} (2008).
\newblock Amenable actions and almost invariant sets.
\newblock {\em Proc. Amer. Math. Soc.} {\bf 136}, 687--697 (electronic).

\bibitem[Kesten (1959)]{MR22:253}
\bibauthor{Kesten, H.} (1959).
\newblock Symmetric random walks on groups.
\newblock {\em Trans. Amer. Math. Soc.} {\bf 92}, 336--354.

\bibitem[K{\l}opotowski, Nadkarni, Sarbadhikari, and Srivastava
  (2002)]{MR1924811}
\bibauthor{K{\l}opotowski, A., Nadkarni, M.G., Sarbadhikari, H., \and{}
  Srivastava, S.M.} (2002).
\newblock Sets with doubleton sections, good sets and ergodic theory.
\newblock {\em Fund. Math.} {\bf 173}, 133--158.

\bibitem[K\"onig (1916)]{Konig:finite}
\bibauthor{K\"onig, D.} (1916).
\newblock {\"U}ber {G}raphen und ihre {A}nwendung auf {D}eterminantentheorie
  und {M}engenlehre.
\newblock {\em Math. Annalen} {\bf 77}, 453--465.

\bibitem[K\"onig (1926)]{Konig:infinite}
\bibauthor{K\"onig, D.} (1926).
\newblock Sur les correspondences multivoques des ensembles.
\newblock {\em Fund. Math.} {\bf 8}, 114--134.

\bibitem[Laczkovich (1988)]{MR89f:28018}
\bibauthor{Laczkovich, M.} (1988).
\newblock Closed sets without measurable matching.
\newblock {\em Proc. Amer. Math. Soc.} {\bf 103}, 894--896.

\bibitem[Lauer and Wormald (2007)]{MR2354714}
\bibauthor{Lauer, J. \and{} Wormald, N.} (2007).
\newblock Large independent sets in regular graphs of large girth.
\newblock {\em J. Combin. Theory Ser. B} {\bf 97}, 999--1009.

\bibitem[Levin, Peres, and Wilmer (2009)]{MR2466937}
\bibauthor{Levin, D.A., Peres, Y., \and{} Wilmer, E.L.} (2009).
\newblock {\em Markov Chains and Mixing Times}.
\newblock American Mathematical Society, Providence, RI.
\newblock With a chapter by James G. Propp and David B. Wilson.

\bibitem[Lov{\'a}sz (1979)]{MR514926}
\bibauthor{Lov{\'a}sz, L.} (1979).
\newblock On the {S}hannon capacity of a graph.
\newblock {\em IEEE Trans. Inform. Theory} {\bf 25}, 1--7.

\bibitem[Ornstein and Weiss (1980)]{MR80j:28031}
\bibauthor{Ornstein, D.S. \and{} Weiss, B.} (1980).
\newblock Ergodic theory of amenable group actions. {I}. {T}he {R}ohlin lemma.
\newblock {\em Bull. Amer. Math. Soc. (N.S.)} {\bf 2}, 161--164.

\bibitem[Tim{\'a}r (2004)]{Timar:pp}
\bibauthor{Tim{\'a}r, {\'A}.} (2004).
\newblock Tree and grid factors for general point processes.
\newblock {\em Electron. Comm. Probab.} {\bf 9}, 53--59 (electronic).

\bibitem[Tim\'ar (2009)]{Timar:IM}
\bibauthor{Tim\'ar, A.} (2009).
\newblock Invariant matchings of exponential tail on coin flips in {$\Z^d$}.
\newblock Preprint.

\bibitem[Tim\'ar (2010)]{Timar:vor}
\bibauthor{Tim\'ar, A.} (2010).
\newblock Invariant colorings of random planar maps.
\newblock {\em Ergodic Theory Dynamical Systems}.
\newblock Available on CJO 22 Apr 2010, doi:10.1017/S0143385709001205.

\bibitem[Wormald (1999)]{MR1725006}
\bibauthor{Wormald, N.C.} (1999).
\newblock Models of random regular graphs.
\newblock In Lamb, J.D. \and{} Preece, D.A., editors, {\em Surveys in
  Combinatorics, 1999}, volume 267 of {\em London Math. Soc. Lecture Note
  Ser.}, pages 239--298. Cambridge Univ. Press, Cambridge.
\newblock Papers from the British Combinatorial Conference held at the
  University of Kent at Canterbury, Canterbury, 1999.

\endreferences
\bye